\documentclass[12pt,reqno]{article}
\usepackage{amssymb,amscd,amsmath,amsthm,color}
\newtheorem{theorem}{Theorem}[section]

\newtheorem{cor}[theorem]{Corollary}
\newtheorem{prop}[theorem]{Proposition}
\newtheorem{conjecture}[theorem]{Conjecture}
\usepackage{enumerate,amssymb}
\theoremstyle{definition}

\newtheorem{question}[theorem]{Question}
\theoremstyle{remark}
\newtheorem{remark}[theorem]{Remark}

\numberwithin{equation}{section}

\def\eps{\varepsilon}

\def\bC{\mathbb{C}}

\def\bM{\mathbb{M}}

\begin{document}
\baselineskip=15pt

\title{ Eigenvalue inequalities for positive block matrices with the inradius of the numerical range.}

\author{ Jean-Christophe Bourin{\footnote{Funded by the ANR Projet (No.\ ANR-19-CE40-0002) and by the French Investissements
 d'Avenir program, project ISITE-BFC (contract ANR-15-IDEX-03).
}}\,  and Eun-Young Lee{\footnote{This research was supported by
Basic Science Research Program through the National Research
Foundation of Korea (NRF) funded by the Ministry of
Education (NRF-2018R1D1A3B07043682)}  }  }

\date{ }

\maketitle

\vskip 10pt\noindent
{\small 
{\bf Abstract.} We prove the operator norm inequality, for a positive matrix partitioned into four blocks in $\bM_n$, 
$$
\left\| \begin{bmatrix} A &X \\
X^* &B\end{bmatrix}\right\|_{\infty} \le \| A+B  \|_{\infty} +\delta(X),
$$
where $\delta(X)$ is the diameter of the largest possible disc in the  numerical range of $X$. This shows that
the inradius   $\varepsilon(X):=\delta(X)/2$ satisfies
$
\varepsilon(X) \ge \| X\|_{\infty} - \| (|X^*|+ |X|)/2\|_{\infty}.
 $ 
Several  eigenvalue inequalities are  derived. In particular, if $X$ is a normal matrix whose spectrum lies in a disc of radius $r$, the third eigenvalue of the full matrix is bounded by the second eigenvalue of the sum of the diagonal block,
$$
\lambda_{3}\left(\begin{bmatrix} A &X \\
X^* &B\end{bmatrix}\right)  \le \lambda_{2}( A+B ) + r.
$$
 We think that $r$ is optimal and we propose a conjecture related to a norm inequality of Hayashi. 

}

\section{Introduction}

Positive matrices partitioned into four blocks play a central role in Matrix Analysis, and in applications, for instance quantum information theory. A lot of 
important theorems deal with these matrices.  Some of these results give comparison between the full matrix and its diagonal blocks, in particular the sum of the diagonal blocks (the partial trace in the quantum terminology). This note focuses on a recent result  of Bourin and Mhana \cite{BM}, involving the numerical range  of the offdiagonal block.   Recall that a symmetric norm $\|\cdot\|$ on $\bM_{2n}$ means a unitarily invariant norm.  It induces a symmetric norm on $\bM_n$ in an obvious way. The Schatten $p$-norms $\|\cdot\|_p$, $1\le p\le \infty$, and the operator norm $(p=\infty)$ are classical examples of symmetric norms. The main result of \cite{BM} reads as follows.

\vskip 5pt
\begin{theorem}\label{thBM} Let $\begin{bmatrix} A &X \\
X^* &B\end{bmatrix} $ be a positive matrix  partitioned into four blocks
in $\bM_n$. Suppose that $W(X)$ has the width $\omega$. Then, for all symmetric norms,
$$
\left\| \begin{bmatrix} A &X \\
X^* &B\end{bmatrix}\right\| \le \| A+B +\omega I \|.
$$
\end{theorem}

\vskip 5pt
Here $I$ stands for the identity matrix, $W(X)$ denotes the numerical range of $X$,  and the width of $W(X)$ is the smallest distance between two parallel straight lines such that the strip between these two lines contains $W(X)$. If $\omega=0$, that is $W(X)$ is a line segment, Theorem 1.1 was first proved by Mhanna \cite{Mha}. Recently \cite{BL2}, Theorem \ref{thBM} has been completed with the reversed inequality
$$
\left\| \begin{bmatrix} A &X \\
X^* &B\end{bmatrix}\right\| \ge \left\|  \begin{bmatrix} \frac{A+B}{2} +d I& 0 \\ 0&\frac{A+B}{2} -d I  \end{bmatrix}\right\|.
$$
where $d:=\min\{ |z|\, :\, z\in W(X)\}$ is the distance from $0$ to $W(X)$. Several applications
were derived.

Some equality cases in Theorem \ref{thBM} occur for the operator norm $\|\cdot\|_{\infty}$ with the following block matrices, where $a,b$ are two arbitrary nonnegative real numbers.
\begin{equation*}\label{ex1}
\begin{bmatrix} \begin{pmatrix} a&0 \\ 0&b \end{pmatrix} &   \begin{pmatrix} 0&a \\ b&0 \end{pmatrix}\\
 \begin{pmatrix} 0&b \\ a&0 \end{pmatrix} &  \begin{pmatrix} b&0 \\ 0&a \end{pmatrix}\end{bmatrix}.
\end{equation*}
This follows from the fact that 
$
W\left(\begin{pmatrix} 0&b \\ a&0 \end{pmatrix} \right)
$
has the width $2\left| |a|-|b|\right|$, a consequence of the classical elliptical range theorem (see \cite{Li} for a short proof).

Though Theorem \ref{thBM} is sharp for the operator norm, a subtle improvement is possible. This is our concern in the next section.
Once again, a geometric feature of $W(X)$ will contribute: its inradius.  Our approach leads to a remarkable list of eigenvalue that cannot be derived from the norm inequalities of Theorem \ref{thBM}.
 The last section is devoted  to some related operator norm inequalities, in particular we will discuss a recent property due to Hayashi \cite{Hay} and propose a conjecture.

\section{Eigenvalue inequalities}

We define the indiameter $\delta(\Lambda)$ of a compact convex set $\Lambda\subset\bC$  as the diameter of the largest possible  disc in $\Lambda$. For matrices $X\in\bM_n$, we shorten
$\delta(W(X))=:\delta(X)$. Recall that the numerical range of a two-by-two matrix is an elliptical disc (or a line segment, or a single point), see Li's short paper \cite{Li} or the famous book \cite{HJ}.

A matrix  $X\in\bM_n$ is identified as an operator on $\bC^n$. If ${\mathcal{S}}$ is a subspace of $\bC^n$, we denote by $X_ {\mathcal{S}}$ the compression of $X$ onto ${\mathcal{S}}$.
We then define the {\it elliptical width} of  $X$ as
$$
\delta_2(X):= \sup_{\dim {\mathcal{S}}=2} \delta(X_ {\mathcal{S}}).
$$
Of course $\delta_2(X)\le \delta(X)\le \omega$ where $\omega$ still denotes the width of $W(X)$. If $X$ is a contraction, then $\delta_2(X)\le 1$, while $\delta(X)$ may be arbitrarily close to 2 (letting $n$ be large enough). We state our main result.

\vskip 5pt
\begin{theorem}\label{thinner} Let $\begin{bmatrix} A &X \\
X^* &B\end{bmatrix} $ be a positive matrix  partitioned into four blocks
in $\bM_n$.  Then, for all $j\in\{0,1,\ldots, n-1\}$,
$$
\lambda_{1+2j}\left(\begin{bmatrix} A &X \\
X^* &B\end{bmatrix}\right)  \le \lambda_{1+j}( A+B ) + \delta_2(X).
$$
\end{theorem}

\vskip 5pt
Here $\lambda_1(S)\ge\cdots\ge \lambda_d(S)$ stand for the eigenvalue of any Hermitian matrix $S\in\bM_d$. If we denote by $\lambda_1^{\uparrow}(S)\le \cdots\le \lambda^{\uparrow}_d(S)$ these eigenvalues arranged in the increasing order, then Theorem \ref{thinner} reads as
$$
\lambda_{2k}^{\uparrow}\left(\begin{bmatrix} A &X \\
X^* &B\end{bmatrix}\right)  \le \lambda_{k}^{\uparrow}( A+B ) + \delta_2(X).
$$
for all $k\in\{1,2,\ldots, n\}$.

The case $j=0$ in Theorem \ref{thinner} improves Theorem \ref{thBM} for the operator norm. We may consider that Theorem \ref{thinner} is trivial for $j=n-1$. Indeed, using the  decomposition \cite[Lemma 3.4]{BL1},
\begin{equation*}
\begin{bmatrix} A &X \\
X^* &B\end{bmatrix} = U
\begin{bmatrix} A &0 \\
0 &0\end{bmatrix} U^* +
V\begin{bmatrix} 0 &0 \\
0 &B\end{bmatrix} V^*,
\end{equation*}
for some unitary matrices $U,\,V\in  \bM_{2n}$, we obtain from Weyl's inequality 
\cite[p.\ 62]{Bh},
\begin{align*}
\lambda_{2n-1}\left(\begin{bmatrix} A &X \\
X^* &B\end{bmatrix}\right)  &\le \lambda_{n} \left(\begin{bmatrix} A &0 \\
0 &0\end{bmatrix}\right)+ \lambda_{n} \left(\begin{bmatrix} 0 &0 \\
0 &B\end{bmatrix}\right) \\
&=  \lambda_{n}(A) + \lambda_n(B) \\ 
&\le \lambda_n(A+B).
\end{align*}

\vskip 5pt
We turn to the proof of the theorem. 

\vskip 5pt
\begin{proof} We first consider the case $j=0$. We may assume that the norm of the block matrix is strictly greater than
the norms of its two diagonal blocks $A$ and $B$, otherwise the statement is trivial. Hence we have two nonzero (column) vectors $h_1, h_2\in\bC^n$ such that $\| h_1\|^2+\| h_2\|^2=1$ and
$$
\lambda_1\left( \begin{bmatrix} A &X \\
X^* &B\end{bmatrix}\right)=
\left\| \begin{bmatrix} A &X \\
X^* &B\end{bmatrix}\right\|_{\infty}= 
\begin{pmatrix}
h^*_1 & h_2^*
\end{pmatrix}
\begin{bmatrix} A &X \\
X^* &B\end{bmatrix}
\begin{pmatrix} h_1 \\ h_2
\end{pmatrix}.
$$
Therefore, denoting by $E_1$ and $E_2$ the rank one projections corresponding to the one dimensional subspaces spanned by $h_1$ and by $h_2$, we have
$$
\left\| \begin{bmatrix} A &X \\
X^* &B\end{bmatrix}\right\|_{\infty}= 
\left\|
\begin{bmatrix} E_1 &0 \\
0 &E_2\end{bmatrix}
\begin{bmatrix} A &X \\
X^* &B\end{bmatrix}
\begin{bmatrix} E_1 &0 \\
0 &E_2\end{bmatrix}
\right\|_{\infty}
$$
Hence, denoting by $F$ a rank two projection such that $E_1\le F$ and $E_2\le F$, we have
\begin{align*}
\left\| \begin{bmatrix} A &X \\
X^* &B\end{bmatrix}\right\|_{\infty}&= 
\left\|
\begin{bmatrix} F &0 \\
0 &F\end{bmatrix}
\begin{bmatrix} A &X \\
X^* &B\end{bmatrix}
\begin{bmatrix} F &0 \\
0 &F\end{bmatrix}
\right\|_{\infty} \\
&=\left\| \begin{bmatrix} FAF &FXF \\
FX^*F &FBF\end{bmatrix}\right\|_{\infty}.
\end{align*}
So, letting ${\mathcal{S}}$  denote the range of $F$, we have
$$
\left\| \begin{bmatrix} A &X \\
X^* &B\end{bmatrix}\right\|_{\infty}
=\left\| \begin{bmatrix} A_{\mathcal{S}} &X_{\mathcal{S}} \\
X^*_{\mathcal{S}} &B_{\mathcal{S}}\end{bmatrix}\right\|_{\infty}.
$$ 
Hence applying Theorem \ref{thBM} for the operator norm, we obtain 
$$
\left\| \begin{bmatrix} A &X \\
X^* &B\end{bmatrix}\right\|_{\infty}\le
\| A_{\mathcal{S}} +B_{\mathcal{S}}  \|_{\infty} + \eps
$$ 
where $\eps$ is the width of $W(X_{\mathcal{S}})$. Since $W(X_{\mathcal{S}})$ is an elliptical disc (as $X_{\mathcal{S}}$ acts on a two-dimensional space), its width equals to its indiameter, hence $\eps\le \delta_2(X)$, and since
$$
\| A_{\mathcal{S}} +B_{\mathcal{S}}  \|_{\infty} =\| (A+B)_{\mathcal{S}} \|_{\infty} \le \| A+B \|_{\infty}=\lambda_1(A+B),
$$
the proof for $j=0$ is complete.

We turn to the general case, $j=1,\ldots, n-1$. By the min-max principle,
\begin{align*}
\lambda_{1+2j}\left(\begin{bmatrix} A &X \\
X^* &B\end{bmatrix}\right)  &\le \inf_{\dim{\mathcal{S}}=n-j}\lambda_{1}\left(\begin{bmatrix} A &X \\
X^* &B\end{bmatrix}_{{\mathcal{S}}\oplus{\mathcal{S}}}\right)  \\
&= \inf_{\dim{\mathcal{S}}=n-j}\lambda_{1}\left(\begin{bmatrix} A_{\mathcal{S}} &X_{\mathcal{S}} \\
X^*_{\mathcal{S}} &B_{\mathcal{S}}\end{bmatrix}\right),
\end{align*}
hence, from the first part of the proof,
\begin{align*}
\lambda_{1+2j}\left(\begin{bmatrix} A &X \\
X^* &B\end{bmatrix}\right)  &\le 
\inf_{\dim{\mathcal{S}}=n-j}\lambda_{1}\left( A_{\mathcal{S}} +
B_{\mathcal{S}}\right) + \delta_2(X) \\
&=\lambda_{1+j}(A+B) + \delta_2(X)
\end{align*}
which is the desired claim. \end{proof}

\vskip 5pt
If $X\in\bM_n$, we denote by ${\mathrm{dist}}(X,\bC I)$ the $\|\cdot\|_{\infty}$-distance from $X$ to  $\bC I$.
Thus, for a  scalar perturbation of a contraction,   $X=\lambda I + C$ for some contraction $C\in\bM_n$ and some $\lambda\in\bC$, we have  ${\mathrm{dist}}(X,\bC I)\le 1$.

\vskip 5pt
\begin{cor}\label{pertu1} Let $\begin{bmatrix} A &X \\
X^* &B\end{bmatrix} $ be a positive matrix  partitioned into four blocks
in $\bM_n$. Then, for all $j\in\{0,1,\ldots, n-1\}$, 
$$
\lambda_{1+2j}\left( \begin{bmatrix} A &X \\
X^* &B\end{bmatrix}\right)  \le \lambda_{1+j} (A+B) + {\mathrm{dist}}(X,\bC I).
$$
\end{cor}

\vskip 5pt
\begin{proof} For any subspace ${\mathcal{S}}\subset C^n$, we have 
$$
{\mathrm{dist}}(X,\bC I)\ge {\mathrm{dist}}(X_{\mathcal{S}},\bC I_{\mathcal{S}}).
$$
If ${\mathcal{S}}$ has dimension 2, then 
$$
{\mathrm{dist}}(X_{\mathcal{S}},\bC I_{\mathcal{S}} )\ge \delta\left(W(X_{\mathcal{S}})\right).
$$ 
Therefore ${\mathrm{dist}}(X,\bC I)\ge \delta_2(X)$ and Theorem \ref{thinner} completes the proof.
\end{proof}

\vskip 5pt
\begin{cor}\label{cor0}  Let $A,B\in\bM_n$. Then, for every $j\ge 0$ such that $1+2j\le n$,
$$
\lambda_{1+2j}\left( A^*A + B^*B \right) \le \lambda_{1+j}\left(AA^*+ BB^*  \right)+\delta_2(AB^*)
$$
\end{cor}

\vskip 5pt
\begin{proof} 
Note that 
$$
\lambda_{1+2j}\left(A^*A + B^*B \right) = \lambda_{1+2j}\left( T^*T\right)  =  \lambda_{1+2j}\left( TT^*\right)  
$$
with $T= \begin{bmatrix} A \\ B \end{bmatrix}$ and
$
TT^*= \begin{bmatrix} AA^* &AB^* \\
BA^* &BB^*\end{bmatrix}
$
so that Theorem \ref{thinner} yields the  desired claim.
\end{proof}

\vskip 5pt
\begin{cor}\label{cornormal1} Let $\begin{bmatrix} A &N \\
N^* &B\end{bmatrix} $ be a positive matrix  partitioned into four blocks
in $\bM_n$. If $N$ is  normal and its spectrum is contained in a disc of radius $r$, then, 
$$
\lambda_{1+2j}\left(\begin{bmatrix} A &N \\
N^* &B\end{bmatrix}\right)  \le \lambda_{1+j}( A+B ) + r.
$$
for all $j=0,1,\ldots,n-1$.
 \end{cor}

\vskip 5pt
\begin{proof}
Corollary \ref{cornormal1} is a special case of corollary \ref{pertu1}, as $N=\lambda I +R$,
where $\lambda$ is the center of the disc of radius $r$ containing the spectrum of $N$, and $\| R\|_{\infty}\le r$.
\end{proof}

\vskip 5pt
\begin{question}\label{quest1} Fix $r>0$ and $\varepsilon>0$. Can we find (with $n$ large enough) a normal matrix $N$  with spectrum in a disc of radius $r$ and  a positive block matrix $\begin{bmatrix} A &N \\
N^* &B\end{bmatrix} $ 
such that
$$
\lambda_{1+2j}\left(\begin{bmatrix} A &N \\
N^* &B\end{bmatrix}\right)  \ge \lambda_{1+j}( A+B ) + r-\varepsilon
$$
for some $j\in\{0,\ldots, n-1\}$ ? Is it true for for $j=0$ ?
\end{question}

\section{Norm inequalities }

Corollary \ref{cornormal1} with $j=0$ reads as follows.

\vskip 5pt
\begin{cor}\label{cornormal} Let $\begin{bmatrix} A &N \\
N^* &B\end{bmatrix} $ be a positive matrix  partitioned into four blocks
in $\bM_n$. If $N$ is  normal and its spectrum is contained in a disc of radius $r$, then, 
$$
\left\| \begin{bmatrix} A &N \\
N^* &B\end{bmatrix}\right\|_{\infty}  \le \| A+B \|_{\infty}  + r.
$$
 \end{cor}

\vskip 5pt
 We do not know wether the constant $r$ is sharp or not (Question \ref{quest1}). If $n=2$, we can replace $r$ by $0$ as the numerical range of $N$ is then a line segment. If $n=3$ there are some simple examples with $N=U$ unitary such that
$$
\left\| \begin{bmatrix} A &U \\
U^* &B\end{bmatrix}\right\|_{\infty}  > \| A+B \|_{\infty}.
$$
See Hayashi's example  in the discussion of \cite[Problem 3]{Hay} and  the interesting study and examples in \cite{GLRT} where we further have $A+B=kI$ for some scalars $k$. The next result is due to Hayashi \cite[Theorem 2.5]{Hay}.

\vskip 5pt
\begin{theorem}\label{thHa} Suppose that $X\in\bM_n$ is invertible with $n$ distinct singular values. If
the inequality 
$$
\left\| \begin{bmatrix} A &X \\
X^* &B\end{bmatrix}\right\|_{\infty}  \le \| A+B \|_{\infty}
$$
holds for all positive block-matrix with $X$ as off-diagonal block, then $X$ is normal.
\end{theorem}

\vskip 5pt
Theorem \ref{thHa} and Theorem \ref{thBM} suggest a natural conjecture. If $W(T)$ is line segment, then $T$ is a so-called  essentially Hermitian matrix.

\vskip 5pt
\begin{conjecture}\label{quest2} Let $X\in\bM_n$. If
the inequality 
$$
\left\| \begin{bmatrix} A &X \\
X^* &B\end{bmatrix}\right\|_{\infty}  \le \| A+B \|_{\infty}
$$
holds for all positive block-matrix with $X$ as off-diagonal block, then $X$ is essentially Hermitian.
\end{conjecture}

\vskip 5pt
 If we replace the operator norm by the Frobenius (or Hilbert-Schmidt) norm $\|\cdot\|_2$ then the following characterization holds.

\vskip 5pt
\begin{prop}\label{prop1} Let $X\in \bM_n$. Then, 
the inequality 
$$
\left\| \begin{bmatrix} A &X \\
X^* &B\end{bmatrix}\right\|_{2}  \le \| A+B \|_{2}
$$
holds for all positive block-matrix with $X$ as off-diagonal block if and only if $X$ is normal.
 \end{prop}

\vskip 5pt
\begin{proof} Suppose that $X$ is normal. To prove the inequality, squaring both side,  it suffices to establish the trace inequality
\begin{equation}\label{tr}
{\mathrm{Tr\,}} X^*X \le {\mathrm{Tr\,}} AB.
\end{equation}
Note that $X=A^{1/2}KB^{1/2}$, for some contraction $K$. Recall that, for all symmetric norms on $\bM_n,$ and any normal matrix $N\in\bM_n$, decomposed as $N=ST$, we have
$\| N\| \le \| TS \|$ (\cite[p.\ 253]{Bh}). Therefore
$$
\| X\| = \| A^{1/2}KB^{1/2}\| \le\| KB^{1/2}A^{1/2}\|.
$$
Squaring this inequality with the Frobenius norm yields the desired  inequality \eqref{tr}.

Suppose that $X$ is nonnormal, and note
that
$
 \begin{bmatrix} |X^*|&X \\
X^* &|X|\end{bmatrix}
$ is positive semidefinite and satisfies 
$$
\left\| \begin{bmatrix} |X^*|&X \\
X^* &|X|\end{bmatrix}
\right\|_2^2= 4 \| |X| \|_2^2
$$
while
$$
\| |X^*| + |X| \|_2^2= 2 \| |X| \|_2^2 + 2 {\mathrm{Tr\,}} |X| |X^*|
$$
In the Hilbert space $(\bM_n, \|\cdot\|_2)$, the assumption
$$
\| |X| \|_2=\| |X^*| \|_2,\quad |X| \neq |X^*|
$$
ensures strict inequality in the Cauchy-Schwarz inequality
$$
{\mathrm{Tr\,}} |X| |X^*| < \| X\|_2^2.
$$
Therefore 
$$
\| |X^*| + |X| \|_2^2 < \left\| \begin{bmatrix} |X^*|&X \\
X^* &|X|\end{bmatrix}
\right\|_2^2
$$
and this completes the proof.
\end{proof}

\vskip 5pt
Proposition \ref{prop1} suggests a question: for which $p\in[1,\infty]$, the schatten $p$-norm inequality
$$
\left\| \begin{bmatrix} A &N \\
N^* &B\end{bmatrix}\right\|_p  \le \| A+B \|_p
$$
holds for any positive partitioned matrices with a normal off-diagonal block $N$ ?

\vskip 5pt
\begin{cor}\label{cor01}  Let $H,K,X\in\bM_n$ be Hermitian. If $X$ is invertible and  $HK$ is a scalar perturbation of a contraction, then,
$$
\left\| XH^2X+ X^{-1}K^2X^{-1} \right\|_{\infty}  \le \left\| HX^2H+ KX^{-2}K  \right\|_{\infty}  + 1.
$$
\end{cor}

\vskip 5pt
\begin{proof} We apply Corollary \ref{cor0} with $j=0$ and $A=HX$, $B=KX^{-1}$, to get
$$
\left\| XH^2X+ X^{-1}K^2X^{-1} \right\|_{\infty}  \le \left\| HX^2H+ KX^{-2}K  \right\|_{\infty}  + \delta_2(HK)
$$
Since $(HK)_{\mathcal{S}}$ is a scalar perturbation of a contraction acting on a space of dimension 2, necessarily $\delta_2(HK)\le 1$.
\end{proof}

\vskip 5pt
For a normal operator,  the numerical range is the convex hull of the spectrum. For a non normal operator $X$,  several lower bounds for the indiameter
of $W(X)$ can be obtained from the left and right modulus $|X^*|$ and $|X|$.

\vskip 5pt
\begin{cor}\label{corcomp} Let $X \in \bM_n$ and let $f(t)$ and $g(t)$ are two nonnegative functions defined on $[0,\infty)$ such that $f(t)g(t)=t^2$ and $f(0)=g(0)=0$. Then,
$$
\delta_2(X) \ge \left\|  f(|X|)+g(|X|) \right\|_{\infty} - \left\|   f(|X^*|) + g(|X|) \right\|_{\infty}.
$$
\end{cor}

\vskip 5pt
\begin{proof} First, observe that we have a function $h(t)$ defined on $[0,\infty)$ such that
\begin{equation}\label{h1}
f(t)=th^2(t^{1/2}), \quad g(t)=th^{-2}(t^{1/2}),
\end{equation}
and $h(t)>0$ for all $t\ge0$ (we may, for instance, set $h(0)=1$). Hence $h(T)$ is invertible for any positive $T$, and from the polar decomposition
$$X=|X^{*}|^{1/2} U |X|^{1/2}$$
with a unitary factor $U$, we infer the factorization
$$
X=|X^{*}|^{1/2}h(|X^{*}|^{1/2}) U |X|^{1/2}h^{-1}(|X|^{1/2}).
$$
Thus $X=AB^*$ where
$
A=|X^{*}|^{1/2}h(|X^{*}|^{1/2}) $  and $B^*=U |X|^{1/2}h^{-1}(|X|^{1/2})$.
Therefore Corollary \ref{cor0} yields
$$
\left\| |X^{*}|h^2(|X^{*}|^{1/2}) +U |X|h^{-2}(|X|^{1/2})U^* \right\|_{\infty} \le
\left\| |X^*|h^2(|X^*|^{1/2}) +|X|h^{-2}(|X|^{1/2}) \right\|_{\infty} 
+\delta_2(X)
$$
Using \eqref{h1} and the fact that $\varphi(|X^*|)=U\varphi(|X|)U^*$ for any function $\varphi(t)$ defined on $[0,\infty)$, the proof is complete.
\end{proof}

\vskip 5pt
The following special case shows that Corollary \ref{corcomp} is rather optimal.

\vskip 5pt
\begin{cor}\label{corsimple} If $X \in \bM_n$ has a numerical range of inradius  $\varepsilon(X)$,  then, for all $a\in\bC$,
$$
\varepsilon(X) \ge \left\|  X-aI\right\|_{\infty} - \left\| \frac{ |X-aI| +|X^*-\overline{a}I|}{2} \right\|_{\infty}.
$$
If $X\in\bM_2$ and $a=\tau$ is the normalized trace of $X$, then this inequality is an equality.
\end{cor}

\vskip 5pt
\begin{proof} Applying Corollary \ref{corcomp} with $X-aI$ and $f(t)=g(t)=t$ yields the inequality. If $X\in\bM_2$, then $X$ is unitarily equivalent to
$$
\begin{pmatrix}
\tau & y \\
x & \tau 
\end{pmatrix}.
$$
So
\begin{align*}
\| X -\tau I\|_{\infty}-\left\| \frac{ |X-\tau I| +|X^*-\overline{\tau}I|}{2} \right\|_{\infty}
&=
\left\|
\begin{pmatrix}
|x| & 0 \\
0 & |y| 
\end{pmatrix}
\right\|_{\infty}
-
\frac{1}{2}\left\|
\begin{pmatrix}
|x| +|y| & 0 \\
0 & |x|+|y| 
\end{pmatrix}
\right\|_{\infty} \\
&=\frac{1}{2} \left| |x|-|y| \right| \\
&=\varepsilon(X)
\end{align*}
establishing the desired equality.
\end{proof}

\vskip 5pt
The special case $f(t)=g(t)=t$ in Corollary \ref{corcomp} seems important, we record it as a proposition:

\vskip 5pt
\begin{prop}\label{thcomp} The elliptical width of the numerical range  of  $X \in \bM_n$ satisfies
$$\delta_2(X)\ge 2\| X\|_{\infty}-\||X|+|X^*|\|_{\infty}. $$
\end{prop}

In particular, the inradius $\varepsilon(X)$ of the numerical range of $X$ satisfies $$\varepsilon(X)\ge \| X\|_{\infty}-\|(|X|+|X^*|)/2\|_{\infty}.$$

\vskip 5pt
\begin{remark}
Our results  still hold for operators on infinite dimensional separable Hilbert space (assuming in Corollary \ref{corcomp} that $f(t)$ and
$g(t)$ are Borel functions).
\end{remark}

  \vskip 15pt\noindent
Eun-Young Lee

\noindent
 Department of mathematics, KNU-Center for Nonlinear
Dynamics,

\noindent
Kyungpook National University,

\noindent
 Daegu 702-701, Korea.

\noindent
 Email: eylee89@knu.ac.kr

\vskip 15pt
\noindent
Jean-Christophe Bourin

\noindent
Laboratoire de math\'ematiques, 

\noindent
Universit\'e de Bourgogne Franche-Comt\'e, 

\noindent
25 000 Besan\c{c}on, France.

\noindent
Email: jcbourin@univ-fcomte.fr

\end{document}